\documentclass[thmsa,12pt]{article}
\usepackage{amsfonts}
\usepackage{amssymb}
\usepackage{amsmath}

\input{tcilatex}
\begin{document}

\begin{center}
{\Large Common Fixed Points of Weakly Commuting Multivalued Mappings\ on
Domain of Sets Endowed with Directed Graph}{\large \medskip }

\textbf{Sergei Silvestrov}$^{1}$ and \textbf{Talat Nazir}$^{1,2}$

{\small $^{(1)}$\textit{Division of Applied Mathematics, School of
Education, Culture and Communication, }}

{\small \textit{M\"{a}lardalen University, 72123 V\"{a}ster\aa s, Sweden.}}

$^{(2)}${\small \textit{Department of Mathematics, COMSATS Institute of
Information Technology,}}

{\small \textit{22060 Abbottabad, Pakistan.}}

E-mail{\small : talat@ciit.net.pk, sergei.silvestrov@mdh.se}\bigskip
\end{center}

\noindent
--------------------------------------------------------------------------------------------

\noindent \textit{Abstract:} \ \ In this paper, the existence of coincidence
points and common fixed points for multivalued mappings satisfying certain
graphic $\psi $-contraction contractive conditions with set-valued domain
endowed with a graph, without appealing to continuity, is established. Some
examples are presented to support the results proved herein. Our results
unify, generalize and extend various results in the existing literature.

\noindent \textbf{---------------------------------------------}


\noindent \textit{2000 \ Mathematics Subject Classification: }\texttt{47H10,
54E50, 54H25}.

\noindent \textbf{---------------------------------------------}\medskip

\section{Introduction and preliminaries}

Order oriented fixed point theory is studied in an environment created by a
class of partially ordered sets with appropriate mappings satisfying certain
order condition like monotonicity, expansivity or order continuity.
Existence of fixed points in partially ordered metric spaces has been
studied by Ran and Reurings \cite{RanReurings}. Recently, many researchers
have obtained fixed point results for single and multivalued mappings
defined on partially ordered metrics spaces (see, e.g., \cite{AminiEmami10,
BegButt13, HarjaniSadar09, NietoLopez05}). Jachymski and Jozwik \cite%
{Jachymski}\textbf{\ }introduced a new approach in metric fixed point theory
by replacing the order structure with a graph structure on a metric space.
In this way, the results proved in ordered metric spaces are generalized
(see also \cite{Jachymski2}\textbf{\ }and the reference therein); in fact,
in 2010, Gwodzdz-lukawska and Jachymski \cite{LJ10}, developed the
Hutchinson-Barnsley theory for finite families of mappings on a metric space
endowed with a directed graph. Abbas and Nazir \cite{AN13} obtained some
fixed point results for power graph contraction pair endowed with a graph.
Bojor \cite{Bojor12} proved fixed point theorems for Reich type contractions
on metric spaces with a graph. For more results in this direction, we refer
to \cite{Aleomraninejad,Alfuraidan, Bojor10, Bojor13, CP12, Nicolae} and
reference mentioned therein.$\smallskip $

Beg and Butt \cite{BegButt} proved the existence of fixed points of
multivalued mapping in metric spaces endowed with a graph $G.$ Recently,
Abbas et al., \cite{AAKN15} obtained fixed points of set valued mappings
satisfying certain graphic contraction conditions with set valued domain
endowed with a graph. Nicolae et al. \cite{Nicolae} established some fixed
points of multivalued generalized contractions in metric spaces endowed with
a graph.

The aim of this paper is to prove some coincidence point and common fixed
point results for discontinuous multivalued graphic $\psi $-contractive
mappings defined on the family of closed and bounded subsets of a metric
space endowed with a graph $G.$ These results extend and strengthen various
comparable results in the existing literature \cite{AAKN15, BegButt,
Bojor10, Jachymski, Jachymski2, Nadler} .

Consistent with Jachymski \cite{Jachymski2}, let $(X,d)$ be a metric space
and $\Delta $ denotes the diagonal of $X\times \ X$. Let $G$ be a directed
graph, such that the set $V(G)$ of its vertices coincides with $X$ and $E(G)$
be the set of edges of the graph which contains all loops, that is, $\Delta
\subseteq E(G)$. Also assume that the graph $G$ has no parallel edges and,
thus, one can identify $G$ with the pair $(V(G),E(G))$.$\smallskip $

\noindent \textbf{Definition 1.1.} \cite{Jachymski2} \ An operator $%
f:X\rightarrow X$ is called a Banach $G$-contraction or simply $G$%
-contraction if

\begin{enumerate}
\item[(a)] $f$ preserves edges of $G;$ for each $x,y\in X$ with $(x,y)\in
E(G),$ we have $(f(x),f(y))\in E(G)$,

\item[(b)] $f$ decreases weights of edges of $G$; there exists $\alpha \in
(0,1)$ such that for all $x,y\in X$ with $(x,y)\in E(G)$, we have $%
d(f(x),f(y))\leq \alpha d(x,y).\smallskip $
\end{enumerate}

If $x$ and $y$ are vertices of $G$, then a path in $G$ from $x$ to $y$ of
length $k\in 
%TCIMACRO{\U{2115} }%
%BeginExpansion
\mathbb{N}
%EndExpansion
$ is a finite sequence $\{x_{n}\}$ ( $n\in \{0,1,2,...,k\}$ ) of vertices
such that $x_{0}=x$, $x_{k}=y$ and $(x_{i-1},x_{i})\in E(G)$ for $i\in
\{1,2,...,k\}$.$\smallskip $

Notice that a graph $G$ is connected if there is a directed path between any
two vertices and it is weakly connected if $\widetilde{G}$ is connected,
where $\widetilde{G}$ denotes the undirected graph obtained from $G$ by
ignoring the direction of edges. Denote by $G^{-1}$ the graph obtained from $%
G$ by reversing the direction of edges. Thus,%
\begin{equation*}
E\left( G^{-1}\right) =\left\{ \left( x,y\right) \in X\times X:\left(
y,x\right) \in E\left( G\right) \right\} .
\end{equation*}%
It is more convenient to treat $\widetilde{G}$ as a directed graph for which
the set of its edges is symmetric, under this convention; we have that%
\begin{equation*}
E(\widetilde{G})=E(G)\cup E(G^{-1}).
\end{equation*}%
In $V(G),$ we define the relation $R$ in the following way:

For $x,y\in V(G),$ we have $xRy$ if and only if, there is a path in $G$ from 
$x$ to $y.$ If $G$ is such that $E(G)$ is symmetric, then for $x\in V(G)$,
the equivalence class $[x]_{\widetilde{G}}$ in $V(G)$ defined by the
relation $R$ is $V(G_{x}).$\smallskip

Recall that if $f:X\rightarrow X$ is an operator, then by $F_{f}$ we denote
the set of all fixed points of $f$. Set%
\begin{equation*}
X_{f}:=\{x\in X:(x,f(x))\in E(G)\}.
\end{equation*}%
Jachymski \cite{Jachymski} used the following property:

(P) : for any sequence $\{x_{n}\}$\ in $X$, if $x_{n}\rightarrow x$ as $%
n\rightarrow \infty $ and $(x_{n},x_{n+1})$ $\in E(G),$ then $(x_{n},x)$ $%
\in E(G).\smallskip $

\noindent \textbf{Theorem 1.2.} \cite{Jachymski} Let $(X,d)$ be a complete
metric space and $G$ a directed graph such that $V(G)=X$ and $f:X\rightarrow
X$ a $G$-contraction. Suppose that $E(G)$ and the triplet $(X,d,G)$ have
property (P). Then the following statements hold:

\begin{enumerate}
\item[(i)] $F_{f}$ $\neq \emptyset $ if and only if $X_{f}\neq \emptyset $;

\item[(ii)] if $X_{f}\neq \emptyset $ and $G$ is weakly connected, then $f$
is a Picard operator, i.e., $F_{f}=\{x^{\ast }\}$ and sequence $%
\{f^{n}(x)\}\rightarrow x^{\ast }$ as $n\rightarrow \infty $, for all $x\in
X $;

\item[(iii)] for any $x\in X_{f}$, $f\mid _{\lbrack x]_{\widetilde{G}}}$ is
a Picard operator;

\item[(iv)] if $X_{f}\subseteq E(G)$, then $f$ is a weakly Picard operator,
i.e., $F_{f}\neq \emptyset $ and, for each $x\in X$, we have sequence $%
\{f^{n}(x)\}\rightarrow x^{\ast }\in F_{f}$ as $n\rightarrow \infty $.
\end{enumerate}

For detailed discussion on Picard operators, we refer to Berinde (\cite%
{Berinde, Berinde5}).$\smallskip $

Let $(X,d)$ be a metric space and $CB(X)$ a class of all nonempty closed and
bounded subsets of $X$. For $A,B\in CB(X)$, let%
\begin{equation*}
H(A,B)=\max \{\sup\limits_{b\in B}d(b,A),\sup\limits_{a\in A}d(a,B)\},
\end{equation*}%
where $d(x,B)=\inf \{d(x,b):b\in B\}$ is the distance of a point $x$ to the
set $B$. The mapping $H$ is said to be the Pompeiu-Hausdorff metric induced
by $d$.\medskip

Throughout this paper, we assume that a directed graph $G$ has no parallel
edge and $G$ is a weighted graph in the sense that each vertex $x$ is
assigned the weight\ $d(x,x)=0$ and each edge $(x,y)$ is assigned the weight 
$d(x,y).$ Since $d$ is a metric on $X,$ the weight assigned to each vertex $%
x $ to vertex $y$ need not be zero and, whenever a zero weight is assigned
to some edge $(x,y),$ it reduces to a loop $(x,x)$ having weight $0.$
Further, in Pompeiu-Hausdorff metric induced by metric $d,$ the
Pompeiu-Hausdorff weight assigned to each $U,V\in CB\left( X\right) $ need
not be zero (that is, $H\left( U,V\right) \neq 0$) and, whenever a zero
Pompeiu-Hausdorff weight is assigned to some $U,V\in CB\left( X\right) ,$
then it reduces to $U=V.$\medskip

\noindent \textbf{Definition 1.3.} \cite{AAKN15}\ Let $A$ and $B$ be two
nonempty subsets of $X$. Then by:

\begin{itemize}
\item[(a)] `there is an edge between $A$ and $B$', we mean there is an edge
between some $a\in A$ and $b\in B$ which we denote by $(A,B)\subset E\left(
G\right) .$

\item[(b)] \ `there is a path between $A$ and $B$', we mean that there is a
path between some $a\in A$ and $b\in B$.
\end{itemize}

In $CB(X),$ we define a relation $R$ in the following way:

For $A,B\in CB(X)$, we have $ARB$ if and only if, there is a path between $A$
and $B$.

We say that the relation $R$ on $CB\left( X\right) $ is transitive if there
is a path between $A$ and $B,$ and there is a path between $B$ and $C,$ then
there is a path between $A$ and $C.\smallskip $

Consider the mapping $T:CB(X)\rightarrow CB(X)$ instead of a mapping $T$
from $X$ to $X$ or from $X$ to $CB(X).$

\noindent For mappings $T:CB\left( X\right) \rightarrow CB\left( X\right) ,$
the set $X_{T}$\ is defined as%
\begin{equation*}
X_{T}:=\{U\in CB\left( X\right) :\left( U,T\left( U\right) \right) \subseteq
E(G)\}.
\end{equation*}

Recently, Abbas et al. \cite{AAKN15} gave the following definition.$%
\smallskip $

\noindent \textbf{Definition 1.4.} \ \ Let $T:CB(X)\rightarrow CB(X)$ a be
multivalued mapping. The mapping $T$ is said to be a graph $\phi $%
-contraction if the following conditions hold:

\begin{itemize}
\item[(i)] There is an edge between $A$ and $B$ implies there is an edge
between $T(A)$ and $T(B)$ for all $A,B\in CB(X)$.

\item[(ii)] There is a path between $A$ and $B$ implies there is a path
between $T(A)$ and $T(B)$ for all $A,B\in CB(X)$.

\item[(iii)] There exists an upper semi-continuous and nondecreasing
function $\phi :\mathbb{%
%TCIMACRO{\U{211d} }%
%BeginExpansion
\mathbb{R}
%EndExpansion
}^{+}\rightarrow \mathbb{%
%TCIMACRO{\U{211d} }%
%BeginExpansion
\mathbb{R}
%EndExpansion
}^{+}$ with $\phi (t)<t$ for each $t>0$ such that there is an edge between $%
A $ and $B$ implies that%
\begin{equation}
H\left( T\left( A\right) ,T\left( B\right) \right) \leq \phi (H(A,B))\text{
for all }A,B\in CB\left( X\right) .  \tag{1.1}
\end{equation}
\end{itemize}

\noindent \textbf{Definition 1.5.} \ \ Let $S,T:CB(X)\rightarrow CB(X)$ be
two multivalued mappings. The set $U\in CB(X)$ is said to be a coincidence
point of $S$ and $T$, if $S\left( U\right) =T\left( U\right) .$ Also, a set $%
A\in CB(X)$ is said to be a fixed point of $S$ if $S(A)=A$. The set of all
coincidence points of $S$ and $T$ is denoted by $CP\left( S,T\right) $ and
the set of all fixed points of $S$ is denoted by $Fix\left( S\right) $%
.\medskip

\noindent \textbf{Definition 1.6.} \ \ Two maps $S,T:CB(X)\rightarrow CB(X)$
are said to be weakly compatible if the commute at their coincidence point.

For more details to the weakly compatible maps, we refer the reader to \cite%
{AJ08, Jungck2, Jungck3}.

\noindent A subset $\Gamma $ of $CB\left( X\right) $ is said to be complete
if for any set $X,Y\in \Gamma ,$ there is an edge between $X$ and $%
Y.\smallskip $

Abbas et al. \cite{AAKN15} used the property $P^{\ast }$ stated as follows:
A graph $G$ is said to have property

\begin{enumerate}
\item[(P$^{\ast }$) :] if for any sequence $\{X_{n}\}$ in $CB(X)$ with $%
X_{n}\rightarrow X$ as $n\rightarrow \infty $, there exists edge between $%
X_{n+1}$ and $X_{n}$ for $n\in 
%TCIMACRO{\U{2115} }%
%BeginExpansion
\mathbb{N}
%EndExpansion
$, implies that there is a subsequence $\{X_{n_{k}}\}$ of $\{X_{n}\}$ with
an edge between $X$ and $X_{n_{k}}$ for $n\in 
%TCIMACRO{\U{2115} }%
%BeginExpansion
\mathbb{N}
%EndExpansion
$.$\smallskip $
\end{enumerate}

\noindent \textbf{Theorem 1.6.} \cite{AAKN15}\ Let $(X,d)$ be a complete
metric space endowed with a directed graph $G$ such that $V(G)=X$ and $%
E(G)\supseteq \Delta .$ If $T:CB\left( X\right) \rightarrow CB\left(
X\right) $ is a graph $\phi $-contraction mapping such that the relation $R$
on $CB\left( X\right) $ is transitive, then following statements hold:

\begin{enumerate}
\item[(a)] if $Fix\left( T\right) $ is complete, then the Pompeiu-Hausdorff
weight assigned to the $U,V\in Fix(T)$ is $0$.

\item[(b)] $X_{T}\neq \emptyset $ provided that $Fix\left( T\right) $ $\neq
\emptyset .$

\item[(c)] If $X_{T}\neq \emptyset $ and the weakly connected graph $G$
satisfies the property (P$^{\ast }$), then $T$ has a fixed point.

\item[(d)] $Fix\left( T\right) $ is complete if and only if $Fix\left(
T\right) $ is a singleton.$\smallskip $
\end{enumerate}

\noindent In the sequel, the letters $%
%TCIMACRO{\U{211d} }%
%BeginExpansion
\mathbb{R}
%EndExpansion
$, $%
%TCIMACRO{\U{211d} }%
%BeginExpansion
\mathbb{R}
%EndExpansion
^{+}$ and $%
%TCIMACRO{\U{2115} }%
%BeginExpansion
\mathbb{N}
%EndExpansion
$ denote the set of all real numbers, the set of all positive real numbers
and the set of all natural numbers, respectively.

\noindent We denote $\Psi $ the set of all functions $\psi :%
%TCIMACRO{\U{211d} }%
%BeginExpansion
\mathbb{R}
%EndExpansion
^{+}\rightarrow 
%TCIMACRO{\U{211d} }%
%BeginExpansion
\mathbb{R}
%EndExpansion
^{+},$ where $\psi $ is nondecreasing function with $\sum_{i=1}^{\infty
}\psi ^{n}(t)$ is convergent. It is easy to show that if $\psi \in \Psi ,$
then $\psi \left( t\right) <t$ for any $t>0$.

We now give the following definition:

\noindent \textbf{Definition 1.7.} \ \ Let $(X,d)$ be a metric space endowed
with a directed graph $G$ such that $V(G)=X$, $E(G)\supseteq \Delta $ and
for every $U\ $in $CB(X),$ $\left( S\left( U\right) ,U\right) \subseteq
E\left( G\right) $ and $\left( U,T\left( U\right) \right) \subseteq E\left(
G\right) $. Let $S,T:CB(X)\rightarrow CB(X)$ be two multivalued mappings.
The pair $(S,T)$ of maps is said to be

\begin{enumerate}
\item[(I)] graph $\psi _{1}$-contraction pair if there exists a $\psi \in
\Psi ,$ there is an edge between $A$ and $B$ such that%
\begin{equation*}
H\left( S\left( A\right) ,S\left( B\right) \right) \leq \psi (M_{1}(A,B))%
\text{ holds,}
\end{equation*}%
where%
\begin{eqnarray*}
M_{1}(A,B) &=&\max \{H(T\left( A\right) ,T\left( B\right) ),H(S\left(
A\right) ,T\left( A\right) ),H(S\left( B\right) ,T\left( B\right) ), \\
&&\dfrac{H(S\left( A\right) ,T\left( B\right) )+H(S\left( B\right) ,T\left(
A\right) )}{2}\}.
\end{eqnarray*}

\item[(II)] graph $\psi _{2}$-contraction pair if there exists a $\psi \in
\Psi ,$ there is an edge between $A$ and $B$ such that%
\begin{equation*}
H\left( S\left( A\right) ,S\left( B\right) \right) \leq \psi (M_{2}(A,B))%
\text{ holds,}
\end{equation*}%
where%
\begin{eqnarray*}
M_{2}(A,B) &=&\alpha H(T\left( A\right) ,T\left( B\right) )+\beta H(S\left(
A\right) ,T\left( A\right) )+\gamma H(S\left( B\right) ,T\left( B\right) ) \\
&&+\delta _{1}H(S\left( A\right) ,T\left( B\right) )+\delta _{2}H(S\left(
B\right) ,T\left( A\right) )
\end{eqnarray*}%
and $\alpha ,\beta ,\gamma ,\delta _{1},$ $\delta _{2}\geq 0,$ $\delta
_{1}\leq \delta _{2}$ with $\alpha +\beta +\gamma +\delta _{1}+\delta
_{2}\leq 1$.
\end{enumerate}

\noindent It is obvious that if a pair $(S,T)$ of multivalued mappings on$\
CB(X)\ $is a graph $\psi _{1}$-contraction or graph $\psi _{2}$-contraction
for graph $G$, then pair $(S,T)\ $is also graph $\psi _{1}$-contraction or
graph $\psi _{2}$-contraction respectively, for the graphs $G^{-1}$, $%
\widetilde{G}$ and $G_{0}$, here the graph $G_{0}$ is defined by $%
E(G_{0})=X\times X$.$\smallskip $

\noindent \textbf{Definition 1.8.} \ \ A metric space $(X,d)$ is called an $%
\varepsilon -$chainable metric space for some $\varepsilon >0$ if for given $%
x,y\in X$, there is $n\in 
%TCIMACRO{\U{2115} }%
%BeginExpansion
\mathbb{N}
%EndExpansion
$ and a sequence $\{x_{n}\}$ such that%
\begin{equation*}
x_{0}=x,\text{ }x_{n}=y\text{ and }d(x_{i-1},x_{i})<\varepsilon \text{ for }%
i=1,...,n.
\end{equation*}

For fixed point result of mappings defined on $\varepsilon -$chainable
metric space, we refer to \cite{BegButt} and references mentioned therein.

We also need of the following lemma of Nadler \cite{Nadler} ( see also, \cite%
{AssadKirk} ).\medskip

\noindent \textbf{Lemma 1.9.} \ \ Let $(X,d)$ be a metric space. If $U,V\in
CB(X)$ with $H(U,V)<\varepsilon $, then for each $u\in U$ there exists an
element $v\in V$ such that $d(u,v)<\varepsilon $.

\section{Common Fixed Points}

\noindent In this section, we obtain coincidence point and common fixed
point results for multivalued selfmaps on $CB(X)$ satisfying graph $\psi $%
-contraction conditions endow with a directed graph.

\noindent \textbf{Theorem 2.1.} \ \ Let $(X,d)$ be a metric space endowed
with a directed graph $G$ such that $V(G)=X$, $E(G)\supseteq \Delta $ and $%
S,T:CB\left( X\right) \rightarrow CB\left( X\right) $ a graph $\psi _{1}$%
-contraction pair such that the range of $T$ contains the range of $S$. Then
the following statements hold:

\begin{enumerate}
\item[(i)] $CP\left( S,T\right) \neq \emptyset $ provided that $G$ is weakly
connected with satisfies the property\ (P$^{\ast }$) and $T\left( X\right) $
is complete subspace of $CB\left( X\right) $.

\item[(ii)] if $CP\left( S,T\right) $ is complete, then the
Pompeiu-Hausdorff weight assigned to the $S\left( U\right) $ and $S\left(
V\right) $ is $0$ for all $U,V\in CP\left( S,T\right) $.

\item[(iii)] if $CP\left( S,T\right) $ is complete and $S$ and $T$ are
weakly compatible, then $Fix\left( S\right) \cap Fix\left( T\right) $ is a
singleton.

\item[(iv)] $Fix\left( S\right) \cap Fix\left( T\right) $ is complete if and
only if $Fix\left( S\right) \cap Fix\left( T\right) $ is a singleton.\medskip
\end{enumerate}

\noindent \textit{Proof. \ \ \ }To prove (i), let $A_{0}\ $be an arbitrary
element in $CB(X).$ Since range of $T$ contains the range of $S$, chosen $%
A_{1}\in CB\left( X\right) $ such that $S\left( A_{0}\right) =T\left(
A_{1}\right) .$ Continuing this process, having chosen $A_{n}$ in $CB\left(
X\right) ,$ we obtain an $A_{n+1}$ in $CB\left( X\right) $ such that $%
S(x_{n})=T(x_{n+1})$ for $n\in 
%TCIMACRO{\U{2115} }%
%BeginExpansion
\mathbb{N}
%EndExpansion
.$ The inclusion $\left( A_{n+1},T\left( A_{n+1}\right) \right) \subseteq
E\left( G\right) $ and $\left( T\left( A_{n+1}\right) ,A_{n}\right) =\left(
S\left( A_{n}\right) ,A_{n}\right) \subseteq E\left( G\right) $ implies that 
$\left( A_{n+1},A_{n}\right) \subseteq E\left( G\right) .$

We may assume that $S\left( A_{n}\right) \neq S\left( A_{n+1}\right) $ for
all $n\in 
%TCIMACRO{\U{2115} }%
%BeginExpansion
\mathbb{N}
%EndExpansion
.$ If not, then $S\left( A_{2k}\right) =S\left( A_{2k+1}\right) $ for some $%
k $, implies $T(A_{2k+1})=S\left( A_{2k+1}\right) ,$ and thus $A_{2k+1}\in
CP\left( S,T\right) .$ Now, since $\left( A_{n+1},A_{n}\right) \subseteq
E\left( G\right) $ for all $n\in 
%TCIMACRO{\U{2115} }%
%BeginExpansion
\mathbb{N}
%EndExpansion
,$ and pair $\left( S,T\right) $ form a graph $\psi _{1}$-contraction, so we
have%
\begin{eqnarray*}
H(T\left( A_{n+1}\right) ,T\left( A_{n+2}\right) ) &=&H(S\left( A_{n}\right)
,S\left( A_{n+1}\right) ) \\
&\leq &\psi \left( M_{1}\left( A_{n},A_{n+1}\right) \right) ,
\end{eqnarray*}%
where%
\begin{eqnarray*}
&&M_{1}\left( A_{n},A_{n+1}\right) \\
&=&\max \{H(T\left( A_{n}\right) ,T\left( A_{n+1}\right) ),H\left( S\left(
A_{n}\right) ,T\left( A_{n}\right) \right) ,H\left( S\left( A_{n+1}\right)
,T\left( A_{n+1}\right) \right) , \\
&&\frac{H\left( S\left( A_{n}\right) ,T\left( A_{n+1}\right) \right)
+H\left( S\left( A_{n+1}\right) ,T\left( A_{n}\right) \right) }{2}\} \\
&=&\max \{H(T\left( A_{n}\right) ,T\left( A_{n+1}\right) ),H\left( T\left(
A_{n+1}\right) ,T\left( A_{n}\right) \right) ,H\left( T\left( A_{n+2}\right)
,T\left( A_{n+1}\right) \right) , \\
&&\frac{H\left( T\left( A_{n+1}\right) ,T\left( A_{n+1}\right) \right)
+H\left( T\left( A_{n+2}\right) ,T\left( A_{n}\right) \right) }{2}\} \\
&\leq &\max \{H\left( T\left( A_{n}\right) ,T\left( A_{n+1}\right) \right)
,H\left( T\left( A_{n+1}\right) ,T\left( A_{n+2}\right) \right) , \\
&&\frac{H\left( T\left( A_{n+2}\right) ,T\left( A_{n+1}\right) \right)
+H\left( T\left( A_{n+1}\right) ,T\left( A_{n}\right) \right) }{2}\} \\
&=&\max \{H\left( T\left( A_{n}\right) ,T\left( A_{n+1}\right) \right)
,H\left( T\left( A_{n+1}\right) ,T\left( A_{n+2}\right) \right) \}.
\end{eqnarray*}%
Thus, we have%
\begin{eqnarray*}
H(T\left( A_{n+1}\right) ,T\left( A_{n+2}\right) &\leq &\psi \left( \max
\{H\left( T\left( A_{n}\right) ,T\left( A_{n+1}\right) \right) ,H\left(
T\left( A_{n+1}\right) ,T\left( A_{n+2}\right) \right) \}\right) \\
&=&\psi (H\left( T\left( A_{n}\right) ,T\left( A_{n+1}\right) \right) )
\end{eqnarray*}%
for all $n\in 
%TCIMACRO{\U{2115} }%
%BeginExpansion
\mathbb{N}
%EndExpansion
.$ Therefore for $i=1,2,...,n$, we have%
\begin{eqnarray*}
H(T\left( A_{i-1}\right) ,T\left( A_{i}\right) ) &\leq &\psi
(H(A_{i-1},A_{i})), \\
H(T\left( A_{i-2}\right) ,T\left( A_{i-1}\right) ) &\leq &\psi
(H(A_{i-2},A_{i-1})), \\
&&\cdot \cdot \cdot , \\
H(T\left( A_{0}\right) ,T\left( A_{1}\right) ) &\leq &\psi (H(A_{0},A_{1})),
\end{eqnarray*}%
and so we obtain%
\begin{equation*}
H(T\left( A_{n}\right) ,T\left( A_{n+1}\right) \leq \psi
^{n}(H(A_{0},T\left( A_{1}\right) ))
\end{equation*}%
for all $n\in 
%TCIMACRO{\U{2115} }%
%BeginExpansion
\mathbb{N}
%EndExpansion
.$ Now for $m,n\in 
%TCIMACRO{\U{2115} }%
%BeginExpansion
\mathbb{N}
%EndExpansion
$ with $m>n\geq 1$, we have%
\begin{eqnarray*}
H\left( T\left( A_{n}\right) ,T\left( A_{m}\right) \right) &\leq &H\left(
T\left( A_{n}\right) ,T\left( A_{n+1}\right) \right) +H\left( T\left(
A_{n+1}\right) ,T\left( A_{n+2}\right) \right) \\
&&+...+H\left( T\left( A_{m-1}\right) ,T\left( A_{m}\right) \right) \\
&\leq &\psi ^{n}(H(A_{0},T\left( A_{1}\right) ))+\psi ^{n+1}(H(A_{0},T\left(
A_{1}\right) )) \\
&&+...+\psi ^{m-1}(H(A_{0},T\left( A_{1}\right) )).
\end{eqnarray*}%
By the convergence of the series $\sum_{i=1}^{\infty }\psi
^{i}(H(A_{0},T\left( A_{1}\right) )),$ we get $H\left( T\left( A_{n}\right)
,T\left( A_{m}\right) \right) \rightarrow 0$ as $n,m\rightarrow \infty $.
Therefore $\{T\left( A_{n}\right) \}$ is a Cauchy sequence in $T\left(
X\right) .$ Since $(T\left( X\right) ,d)$ is complete in $CB\left( X\right) $%
, we have $T\left( A_{n}\right) \rightarrow V$ as $n\rightarrow \infty $\
for some $V\in CB\left( X\right) .$ Also, we can find $U$ in $CB\left(
X\right) $ such that $T(U)=V.$

We claim that $S(U)=T(U).$ If not, then since $\left( T\left( A_{n+1}\right)
,T\left( A_{n}\right) \right) \subseteq E\left( G\right) $ so by property (P$%
^{\ast }$), there exists a subsequence $\{T\left( A_{n_{k}+1}\right) \}$ of $%
\{T\left( A_{n+1}\right) \}$ such that $\left( T\left( U\right) ,T\left(
A_{n_{k}+1}\right) \right) \subseteq E\left( G\right) $ for every $n\in 
%TCIMACRO{\U{2115} }%
%BeginExpansion
\mathbb{N}
%EndExpansion
$. As $\left( U,T\left( U\right) \right) \subseteq E\left( G\right) $ and $%
\left( T\left( A_{n_{k}+1}\right) ,A_{n_{k}}\right) =\left( S\left(
A_{n_{k}}\right) ,A_{n_{k}}\right) \subseteq E\left( G\right) $ implies that 
$\left( U,A_{n_{k}}\right) \subseteq E\left( G\right) .$ Now%
\begin{equation}
H(S\left( U\right) ,T\left( A_{n_{k}+1}\right) )=H(S\left( U\right) ,S\left(
A_{n_{k}}\right) )\leq \psi \left( M_{1}\left( U,A_{n_{k}}\right) \right) , 
\tag{(1)}
\end{equation}%
where%
\begin{eqnarray*}
M_{1}\left( U,A_{n_{k}}\right) &=&\max \{H(T\left( U\right) ,T\left(
A_{n_{k}}\right) ),H(S\left( U\right) ,T\left( U\right) ),H(S\left(
A_{n_{k}}\right) ,T\left( A_{n_{k}}\right) ), \\
&&\frac{H(S\left( U\right) ,T\left( A_{n_{k}}\right) )+H(S\left(
A_{n_{k}}\right) ,T\left( U\right) )}{2}\} \\
&=&\max \{H(T\left( U\right) ,T\left( A_{n_{k}}\right) ),H(S\left( U\right)
,T\left( U\right) ),H(T\left( A_{n_{k}+1}\right) ,T\left( A_{n_{k}}\right) ),
\\
&&\frac{H(S\left( U\right) ,T\left( A_{n_{k}}\right) )+H(T\left(
A_{n_{k}+1}\right) ,T\left( U\right) )}{2}\}.
\end{eqnarray*}%
Now we consider the following cases:

\noindent If $M_{1}\left( U,A_{n_{k}}\right) =H(T\left( U\right) ,T\left(
A_{n_{k}}\right) ),$ then on taking limit as $k\rightarrow \infty $ in (1),
we have%
\begin{equation*}
H(S\left( U\right) ,T\left( U\right) )\leq \psi \left( H\left( T\left(
U\right) ,T\left( U\right) \right) \right) ,
\end{equation*}%
a contradiction.

\noindent When $M_{1}\left( U,A_{n_{k}}\right) =H(S\left( U\right) ,T\left(
U\right) ),$ then%
\begin{equation*}
H(S\left( U\right) ,T\left( U\right) )\leq \psi \left( H\left( S\left(
U\right) ,T\left( U\right) \right) \right) ,
\end{equation*}%
gives a contradiction.

\noindent In case $M_{1}\left( U,A_{n_{k}}\right) =H(T\left(
A_{n_{k}+1}\right) ,T\left( A_{n_{k}}\right) ),$ then on taking limit as $%
k\rightarrow \infty $ in (1), we get%
\begin{equation*}
H(S\left( U\right) ,T\left( U\right) )\leq \psi \left( H\left( T\left(
U\right) ,T\left( U\right) \right) \right) ,
\end{equation*}%
a contradiction.

\noindent Finally, if $M_{1}\left( U,A_{n_{k}}\right) =\dfrac{H(S\left(
U\right) ,T\left( A_{n_{k}}\right) )+H(T\left( A_{n_{k}+1}\right) ,T\left(
U\right) )}{2},$ then on taking limit as $k\rightarrow \infty ,$ we have%
\begin{eqnarray*}
H(S\left( U\right) ,T\left( U\right) ) &\leq &\psi (\frac{H\left( S\left(
U\right) ,T\left( U\right) \right) +H\left( T\left( U\right) ,T\left(
U\right) \right) }{2}) \\
&=&\psi (\frac{H\left( S\left( U\right) ,T\left( U\right) \right) }{2}),
\end{eqnarray*}%
a contradiction.

\noindent Hence $S\left( U\right) =T\left( U\right) ,$ that is, $U\in
CP\left( S,T\right) $.

\noindent To prove (ii), suppose that $CP\left( S,T\right) $ is complete set
in $G$. Let $U,V\in CP\left( S,T\right) $ and suppose that the
Pompeiu-Hausdorff weight assign to the $S\left( U\right) $ and $S\left(
V\right) $ is not zero. Since pair $(S,T)$ is a graph $\psi _{1}$%
-contraction, we obtain that%
\begin{eqnarray*}
&&H(S\left( U\right) ,S\left( V\right) ) \\
&\leq &\psi (M_{1}(U,V)) \\
&\leq &\psi (\max \{H(T\left( U\right) ,T\left( V\right) ),H(S\left(
U\right) ,T\left( U\right) ),H(S\left( V\right) ,T\left( V\right) ), \\
&&\frac{H(S\left( U\right) ,T\left( V\right) )+H(S\left( V\right) ,T\left(
U\right) )}{2}\}) \\
&=&\psi (\max \{H\left( S\left( U\right) ,S\left( V\right) \right)
,H(S\left( U\right) ,S\left( U\right) ),H(S\left( V\right) ,T\left( V\right)
), \\
&&\frac{H(S\left( U\right) ,S\left( V\right) )+H(S\left( V\right) ,S\left(
U\right) )}{2}\}) \\
&=&\psi \left( H(S\left( U\right) ,S\left( V\right) )\right) ,
\end{eqnarray*}%
a contradiction as $\psi \left( t\right) <t$ for all $t>0$. Hence (ii) is
proved.

\noindent To prove (iii),\ suppose the set $CP\left( S,T\right) \ $is weakly
compatible.\ First we are to show that $Fix\left( T\right) \cap Fix(S)\ $is
nonempty. Let $W=S\left( U\right) =T\left( U\right) ,$ then we have $T\left(
W\right) =TS\left( U\right) =ST\left( U\right) =S\left( W\right) ,$ which
shows that $W\in CP\left( S,T\right) .$ Thus the Pompeiu-Hausdorff weight
assign to the $S\left( U\right) $ and $S\left( W\right) $ is zero$\ $(by
ii). Hence $W=S\left( W\right) =T\left( W\right) ,$ that is, $W\in Fix\left(
S\right) \cap Fix\left( T\right) .$ Since $CP\left( S,T\right) $ is
singleton set, implies $Fix\left( S\right) \cap Fix\left( T\right) $ is
singleton.

\noindent Finally to prove (iv),\ suppose the set $Fix\left( S\right) \cap
Fix\left( T\right) \ $is complete.\ We are to show that $Fix\left( T\right)
\cap Fix(S)\ $is singleton. Assume on contrary that there exist $U$,$V\in
CB\left( X\right) $ such that $U,V\in Fix\left( S\right) \cap Fix\left(
T\right) $ and $U\neq V$. By completeness of $Fix\left( S\right) \cap
Fix\left( T\right) $,\ there exists an edge between $U$ and $V.$ As pair $%
(S,T)$ is a graph $\psi _{1}$-contraction, so we have%
\begin{eqnarray*}
H(U,V) &=&H(S\left( U\right) ,S\left( V\right) ) \\
&\leq &\psi (M_{1}(U,V)) \\
&=&\psi (\max \{H(T\left( U\right) ,T\left( V\right) ),H(S\left( U\right)
,T\left( U\right) ),H(S\left( V\right) ,T\left( V\right) ), \\
&&\frac{H(S\left( U\right) ,T\left( V\right) )+H(S\left( V\right) ,T\left(
U\right) )}{2}\}) \\
&=&\psi (\max \{H(U,V),H(U,U),H(V,V), \\
&&\frac{H(U,V)+H(V,U)}{2}\}) \\
&=&\psi \left( H\left( U,V\right) \right) ,
\end{eqnarray*}%
a contradiction.\ Hence $U=V$. Conversely, if $Fix(S)\cap Fix(T)$ is
singleton, then since $E(G)\supseteq \Delta ,$ so it is obvious that $%
F(S)\cap F(T)$ is complete set. $\square $\medskip

\noindent \textbf{Example 2.2.} \ \ Let $X=\{1,2,...,n\}=V\left( G\right) ,$ 
$n>2$\ and $E\left( G\right) =\{\left( i,j\right) \in X\times X:$ $i\leq
j\}. $ Let $V\left( G\right) $ be endowed with metric $d:X\times
X\rightarrow 
%TCIMACRO{\U{211d} }%
%BeginExpansion
\mathbb{R}
%EndExpansion
^{+}$ defined by%
\begin{equation*}
d\left( x,y\right) =\left\{ 
\begin{array}{cc}
0 & \text{if }x=y, \\ 
&  \\ 
\dfrac{1}{n} & \text{if }x\in \{1,2\}\text{ with }x\neq y, \\ 
&  \\ 
\dfrac{n}{n+1} & \text{otherwise.}%
\end{array}%
\right.
\end{equation*}%
Furthermore, the Pompeiu-Hausdorff metric is given by%
\begin{equation*}
H(A,B)=\left\{ 
\begin{array}{ll}
\text{ \ }\dfrac{1}{n} & \text{if }A,B\subseteq \{1,2\}\text{ with }A\neq B,
\\ 
&  \\ 
\dfrac{n}{n+1} & \text{if }A\text{ or }B\text{ (or both)}\nsubseteq \{1,2\}%
\text{ with }A\neq B, \\ 
&  \\ 
\text{ \ \ }0 & \text{if }A=B.%
\end{array}%
\right.
\end{equation*}%
The Pompeiu-Hausdorff weights (for $n=4$) assigned to $A,B\in CB\left(
X\right) $ are shown in the Figure.%
\begin{equation*}
\FRAME{itbpF}{4.7806in}{3.4652in}{0in}{}{}{Figure}{\special{language
"Scientific Word";type "GRAPHIC";maintain-aspect-ratio TRUE;display
"USEDEF";valid_file "T";width 4.7806in;height 3.4652in;depth
0in;original-width 8.4056in;original-height 6.0806in;cropleft "0";croptop
"1";cropright "1";cropbottom "0";tempfilename
'NQDDHN00.wmf';tempfile-properties "XPR";}}
\end{equation*}

\noindent Define $S,T:CB\left( X\right) \rightarrow CB(X)$ as follows:%
\begin{eqnarray*}
S(U) &=&\left\{ 
\begin{array}{cc}
\{1\}, & \text{if }U\subseteq \{1,2\}, \\ 
\{1,2\}, & \text{if }U\varsubsetneq \{1,2\}%
\end{array}%
\right. \\
T(U) &=&\left\{ 
\begin{array}{cl}
\{1\}, & \text{if }U=\{1\}, \\ 
\{1,2,3\}, & \text{if }U\subseteq \{2,3\}. \\ 
\{1,2,...,n\}, & \text{otherwise.}%
\end{array}%
\right.
\end{eqnarray*}%
Note that, for all $V\in CB(X),$ $(V,S\left( V\right) )\subseteq E\left(
G\right) $ and $(V,T\left( V\right) )\subseteq E\left( G\right) $.

\noindent Let $\psi :%
%TCIMACRO{\U{211d} }%
%BeginExpansion
\mathbb{R}
%EndExpansion
_{+}\rightarrow 
%TCIMACRO{\U{211d} }%
%BeginExpansion
\mathbb{R}
%EndExpansion
_{+}$ be defined by%
\begin{equation*}
\psi \left( \alpha \right) =\left\{ 
\begin{array}{cc}
\dfrac{1}{2}\alpha ^{2} & 0\leq \alpha <\frac{1}{2} \\ 
&  \\ 
\dfrac{\alpha }{\alpha +1}, & \frac{1}{2}\leq \alpha .%
\end{array}%
\right.
\end{equation*}%
It is easy to verify that $\psi \in \Psi .$\ Now for all $A,B\in CB\left(
X\right) $ with $S\left( A\right) \neq S\left( B\right) ,$ we consider the
following cases:

\begin{description}
\item[(i)] If $A\subseteq \left\{ 1,2\right\} $ and $B=\left\{ 3\right\} $
with $\left( A,B\right) \subseteq E\left( G\right) ,$ then we have%
\begin{eqnarray*}
H(S\left( A\right) ,S\left( B\right) ) &=&H\left( \{1\},\{1,2\}\right) \\
&=&\frac{1}{n} \\
&<&\frac{n}{2n+1} \\
&=&\psi \left( \frac{n}{n+1}\right) \\
&=&\psi \left( H\left( \{1,2\},\{1,2,3\}\right) \right) \\
&=&\psi \left( H(S\left( B\right) ,T\left( B\right) )\right) \leq \psi
(M_{1}(A,B)).
\end{eqnarray*}

\item[(ii)] When $A\subseteq \left\{ 1,2\right\} $ and $B\varsubsetneq
\left\{ 1,2,3\right\} $ with $\left( A,B\right) \subseteq E\left( G\right) ,$
implies that%
\begin{eqnarray*}
H(S\left( A\right) ,S\left( B\right) ) &=&H\left( \{1\},\{1,2\}\right) \\
&=&\frac{1}{n} \\
&<&\frac{n}{2n+1} \\
&=&\psi \left( \frac{n}{n+1}\right) \\
&=&\psi \left( H\left( \{1,2\},\{1,2,...,n\}\right) \right) \\
&=&\psi \left( H(S\left( B\right) ,T\left( B\right) )\right) \leq \psi
(M_{1}(A,B)).
\end{eqnarray*}

\item[(iii)] In case $A=\{3\}$ and $B\subseteq \left\{ 1,2\right\} $ and
with $\left( A,B\right) \subseteq E\left( G\right) ,$ we have%
\begin{eqnarray*}
H(S\left( A\right) ,S\left( B\right) ) &=&H\left( \{1,2\},\{1\}\right) \\
&=&\frac{1}{n} \\
&<&\frac{n}{2n+1} \\
&=&\psi \left( \frac{n}{n+1}\right) \\
&=&\psi \left( H\left( \{1,2\},\{1,2,3\}\right) \right) \\
&=&\psi \left( H(S\left( A\right) ,T\left( A\right) )\right) \leq \psi
(M_{1}(A,B)).
\end{eqnarray*}

\item[(iv)] When $A\varsubsetneq \left\{ 1,2,3\right\} $ and $B\subseteq
\left\{ 1,2\right\} \ $with $\left( A,B\right) \subseteq E\left( G\right) ,$
implies that%
\begin{eqnarray*}
H(S\left( A\right) ,S\left( B\right) ) &=&H\left( \{1,2\},\{1\}\right) \\
&=&\frac{1}{n} \\
&<&\frac{n}{2n+1} \\
&=&\psi \left( \frac{n}{n+1}\right) \\
&=&\psi \left( H\left( \{1,2\},\{1,2,...,n\}\right) \right) \\
&=&\psi \left( H(S\left( A\right) ,T\left( A\right) )\right) \leq \psi
(M_{1}(A,B)).
\end{eqnarray*}
\end{description}

\noindent Hence pair $\left( S,T\right) $\ is graph $\psi _{1}$-contraction.
Thus all the conditions of Theorem 1 are satisfied. Moreover, $\{1\}$ is the
common fixed point of $S$ and $T$, and $Fix\left( S\right) \cap Fix\left(
T\right) $ is complete. $\square $\medskip

\noindent In the next example we show that it is not necessary the given
graph $\left( V\left( G\right) ,E\left( G\right) \right) $ will always be
complete graph.

\noindent \textbf{Example 2.3.} \ \ Let $X=\{1,2,...,n\}=V\left( G\right) ,$ 
$n>2$\ and%
\begin{eqnarray*}
E\left( G\right) &=&\{\left( 1,1\right) ,(2,2),...,(n,n), \\
&&\left( 1,2\right) ,...,(1,n)\}.
\end{eqnarray*}%
On $V\left( G\right) ,$ the metric $d:X\times X\rightarrow 
%TCIMACRO{\U{211d} }%
%BeginExpansion
\mathbb{R}
%EndExpansion
^{+}$ and Pompeiu-Hausdorff metric $H:CB\left( X\right) \rightarrow 
%TCIMACRO{\U{211d} }%
%BeginExpansion
\mathbb{R}
%EndExpansion
^{+}$ are defined as in Example 2.2. The Pompeiu-Hausdorff weights (for $n=4$%
) assigned to $A,B\in CB\left( X\right) $ are shown in the Figure.%
\begin{equation*}
\FRAME{itbpF}{4.6773in}{3.4669in}{0in}{}{}{Figure}{\special{language
"Scientific Word";type "GRAPHIC";maintain-aspect-ratio TRUE;display
"USEDEF";valid_file "T";width 4.6773in;height 3.4669in;depth
0in;original-width 8.3022in;original-height 6.1421in;cropleft "0";croptop
"1";cropright "1";cropbottom "0";tempfilename
'NQDDHN01.wmf';tempfile-properties "XPR";}}
\end{equation*}

\noindent Define $S,T:CB\left( X\right) \rightarrow CB(X)$ as follows:%
\begin{eqnarray*}
S(U) &=&\left\{ 
\begin{array}{cc}
\{1\}, & \text{if }U=\{1\}, \\ 
\{1,2\}, & \text{if }U\neq \{1\}%
\end{array}%
\right. \\
T(U) &=&\left\{ 
\begin{array}{cc}
\{1\}, & \text{if }U=\{1\}, \\ 
\{1,...,n\}, & \text{if }U\neq \{1\}.%
\end{array}%
\right.
\end{eqnarray*}%
Note that, $(S\left( A\right) ,A)\subseteq E\left( G\right) $ and $%
(A,T\left( A\right) )\subseteq E\left( G\right) $ for all $A\in CB(X)$.

\noindent Take $\psi \left( \alpha \right) =\left\{ 
\begin{array}{cc}
\frac{1}{8}t, & t\in \lbrack 0,\frac{1}{4}] \\ 
&  \\ 
\frac{t+1}{t+2}, & t\geq \frac{1}{4}.%
\end{array}%
\right. $. Note that $\psi \in \Psi .$

\noindent For all $A,B\in CB\left( X\right) $ with $S\left( A\right) \neq
S\left( B\right) ,$ we consider the following cases:

\begin{description}
\item[(I)] If $A=\left\{ 1\right\} $ and $B\neq \left\{ 1\right\} ,$ then we
have%
\begin{eqnarray*}
H(S\left( A\right) ,S\left( B\right) ) &=&\frac{1}{n} \\
&<&\frac{2n+1}{3n+1} \\
&=&\psi \left( \frac{n}{n+1}\right) \\
&=&\psi \left( H(S\left( B\right) ,T\left( B\right) )\right) \leq \psi
(M_{1}(A,B)).
\end{eqnarray*}

\item[(II)] If $A\neq \left\{ 1\right\} $ and $B=\left\{ 1\right\} ,$ then
we have%
\begin{eqnarray*}
H(S\left( A\right) ,S\left( B\right) ) &=&\frac{1}{n} \\
&<&\frac{2n+1}{3n+1} \\
&=&\psi \left( \frac{n}{n+1}\right) \\
&=&\psi \left( H(S\left( A\right) ,T\left( A\right) )\right) \leq \psi
\left( M_{1}\left( A,B\right) \right) .
\end{eqnarray*}
\end{description}

\noindent Hence pair $\left( S,T\right) $\ is graph $\psi _{1}$-contraction.
Thus all the conditions of Theorem 1 are satisfied. Moreover, $S$ and $T$
have a common fixed point and $Fix\left( S\right) \cap Fix\left( T\right) $
is complete in $CB\left( X\right) $. $\square $\medskip

\noindent \textbf{Theorem 2.4.} \ \ Let $(X,d)$ be a $\varepsilon -$%
chainable complete metric space for some $\varepsilon >0$ and $S,T:CB\left(
X\right) \rightarrow CB(X)$ be multivalued mappings. Suppose that for all $%
A,B\in CB\left( X\right) ,$%
\begin{equation*}
0<H\left( S\left( A\right) ,S\left( B\right) \right) <\varepsilon
\end{equation*}%
and there exists a $\psi \in \Psi $ such%
\begin{equation*}
H\left( S\left( A\right) ,S\left( B\right) \right) \leq \psi (M_{1}(A,B)),
\end{equation*}%
hold where%
\begin{eqnarray*}
M_{1}(A,B) &=&\max \{H(T\left( A\right) ,T\left( B\right) ),H(S\left(
A\right) ,T\left( A\right) ),H(S\left( B\right) ,T\left( B\right) ), \\
&&\dfrac{H(S\left( A\right) ,T\left( B\right) )+H(S\left( B\right) ,T\left(
A\right) )}{2}\}.
\end{eqnarray*}%
Then $S$ and $T$ have a common fixed point provided that $S$ and $T$ are
weakly compatible.

\noindent \textit{Proof.} \ \ \ By Lemma 1.9,\ from $H\left( A,B\right)
<\varepsilon ,$ we have for each $a\in A,$ an element $b\in B$ such that $%
d(a,b)<\varepsilon $. Consider the graph $G$ as $V(G)=X$ and%
\begin{equation*}
E(G)=\{(a,b)\in X\times X:0<d(a,b)<\varepsilon \}.
\end{equation*}%
Then the $\varepsilon -$chainability of $(X,d)$ implies that $G$ is
connected. For $(A,B)\subset E(G)$, we have from the hypothesis%
\begin{equation*}
H\left( S\left( A\right) ,S\left( B\right) \right) \leq \psi (M_{1}(A,B)),
\end{equation*}%
\begin{eqnarray*}
\text{where }M_{1}(A,B) &=&\max \{H(S\left( A\right) ,T\left( B\right)
),H(S\left( A\right) ,T\left( A\right) ),H(S\left( B\right) ,T\left(
B\right) ), \\
&&\dfrac{H(S\left( A\right) ,T\left( B\right) )+H(S\left( B\right) ,T\left(
A\right) )}{2}\}
\end{eqnarray*}%
implies that pair $(S,T)\ $is graph $\psi _{1}-$contraction.

\noindent Also, $G$ has property (P$^{\ast }$). Indeed, if $\{X_{n}\}$ in $%
CB(X)$ with $X_{n}\rightarrow X$ as $n\rightarrow \infty $ and $\left(
X_{n},X_{n+1}\right) \subset E\left( G\right) $ for $n\in 
%TCIMACRO{\U{2115} }%
%BeginExpansion
\mathbb{N}
%EndExpansion
$, implies that there is a subsequence $\{X_{n_{k}}\}$ of $\{X_{n}\}$ such
that $\left( X_{n_{k}},X\right) \subset E\left( G\right) $ for $n\in 
%TCIMACRO{\U{2115} }%
%BeginExpansion
\mathbb{N}
%EndExpansion
$. So by Theorem 2.1 (iii), $S$ and $T$\ have a common fixed point. $\square 
$\medskip

\noindent \textbf{Corollary 2.5.} \ \ Let $(X,d)$ be a complete metric space
endowed with a directed graph $G$ such that $V(G)=X$ and $E(G)\supseteq
\Delta .$ Suppose that the mapping $S:CB\left( X\right) \rightarrow CB\left(
X\right) $ satisfies the following:

\begin{itemize}
\item[(a)] for every $V\ $in $CB(X),$ $\left( V,S\left( V\right) \right)
\subset E\left( G\right) $.

\item[(b)] There exists $\psi \in \Psi $ such that there is an edge between $%
A$ and $B$ implies that%
\begin{equation*}
H(S\left( A\right) ,S\left( B\right) )\leq \psi (M_{1}(A,B)),
\end{equation*}%
where%
\begin{eqnarray*}
M_{1}(A,B) &=&\max \{H(A,B),H(A,S\left( A\right) ),H(B,S\left( B\right) ), \\
&&\frac{H(A,S\left( B\right) ),H(B,S\left( A\right) )}{2}\}).
\end{eqnarray*}%
Then following statements hold:
\end{itemize}

\begin{enumerate}
\item[(i)] if $Fix\left( S\right) $ is complete, then the Pompeiu-Hausdorff
weight assigned to the $U,V\in Fix\left( S\right) $ is $0$.

\item[(ii)] If the weakly connected graph $G$ satisfies the property (P$%
^{\ast }$), then $S$ has a fixed point.

\item[(iii)] $Fix\left( S\right) \ $is complete if and only if $Fix\left(
S\right) \ $is a singleton.
\end{enumerate}

\noindent \textit{Proof.} \ \ \ Take $T=I$ (identity map) in (1.2), then
Corollary 2.5 follows from Theorem 2.1. $\square $\medskip

\noindent \textbf{Theorem 2.6.} \ \ Let $(X,d)$ be a metric space endowed
with a directed graph $G$ such that $V(G)=X$, $E(G)\supseteq \Delta $ and $%
S,T:CB\left( X\right) \rightarrow CB\left( X\right) $ a graph $\psi _{2}$%
-contraction pair such that the range of $T$ contains the range of $S$. Then
the following statements hold:

\begin{enumerate}
\item[(i)] $CP\left( S,T\right) \neq \emptyset $ provided that $G$ is weakly
connected with satisfies the property\ (P$^{\ast }$) and $T\left( X\right) $
is complete subspace of $CB\left( X\right) $.

\item[(ii)] if $CP\left( S,T\right) $ is complete, then the
Pompeiu-Hausdorff weight assigned to the $S\left( U\right) $ and $S\left(
V\right) $ is $0$ for all $U,V\in CP\left( S,T\right) $.

\item[(iii)] if $CP\left( S,T\right) $ is complete and $S$ and $T$ are
weakly compatible, then $Fix\left( S\right) \cap Fix\left( T\right) $ is a
singleton.

\item[(iv)] $Fix\left( S\right) \cap Fix\left( T\right) $ is complete if and
only if $Fix\left( S\right) \cap Fix\left( T\right) $ is a singleton.\medskip
\end{enumerate}

\noindent \textit{Proof. \ \ \ }To prove (i), let $A_{0}\ $be an arbitrary
element in $CB(X).$ Since range of $T$ contains the range of $S$, chosen $%
A_{1}\in CB\left( X\right) $ such that $S\left( A_{0}\right) =T\left(
A_{1}\right) .$ Continuing this process, having chosen $A_{n}$ in $CB\left(
X\right) ,$ we obtain an $A_{n+1}$ in $CB\left( X\right) $ such that $%
S(x_{n})=T(x_{n+1})$ for $n\in 
%TCIMACRO{\U{2115} }%
%BeginExpansion
\mathbb{N}
%EndExpansion
.$ The inclusion $\left( A_{n+1},T\left( A_{n+1}\right) \right) \subseteq
E\left( G\right) $ and $\left( T\left( A_{n+1}\right) ,A_{n}\right) =\left(
S\left( A_{n}\right) ,A_{n}\right) \subseteq E\left( G\right) $ implies that 
$\left( A_{n+1},A_{n}\right) \subseteq E\left( G\right) .$

We may assume that $S\left( A_{n}\right) \neq S\left( A_{n+1}\right) $ for
all $n\in 
%TCIMACRO{\U{2115} }%
%BeginExpansion
\mathbb{N}
%EndExpansion
.$ If not, then $S\left( A_{2k}\right) =S\left( A_{2k+1}\right) $ for some $%
k $, implies $T(A_{2k+1})=S\left( A_{2k+1}\right) ,$ and thus $A_{2k+1}\in
CP\left( S,T\right) .$ Now, since $\left( A_{n+1},A_{n}\right) \subseteq
E\left( G\right) $ for all $n\in 
%TCIMACRO{\U{2115} }%
%BeginExpansion
\mathbb{N}
%EndExpansion
,$ and pair $\left( S,T\right) $ form a graph $\psi _{2}$-contraction, so we
have%
\begin{eqnarray*}
H(T\left( A_{n+1}\right) ,T\left( A_{n+2}\right) ) &=&H(S\left( A_{n}\right)
,S\left( A_{n+1}\right) ) \\
&\leq &\psi \left( M_{2}\left( A_{n},A_{n+1}\right) \right) ,
\end{eqnarray*}%
where%
\begin{eqnarray*}
&&M_{2}\left( A_{n},A_{n+1}\right) \\
&=&\alpha H(T\left( A_{n}\right) ,T\left( A_{n+1}\right) )+\beta H\left(
S\left( A_{n}\right) ,T\left( A_{n}\right) \right) +\gamma H\left( S\left(
A_{n+1}\right) ,T\left( A_{n+1}\right) \right) \\
&&\delta _{1}H\left( S\left( A_{n}\right) ,T\left( A_{n+1}\right) \right)
+\delta _{2}H\left( S\left( A_{n+1}\right) ,T\left( A_{n}\right) \right) \\
&=&\alpha H(T\left( A_{n}\right) ,T\left( A_{n+1}\right) )+\beta H\left(
T\left( A_{n+1}\right) ,T\left( A_{n}\right) \right) +\gamma H\left( T\left(
A_{n+2}\right) ,T\left( A_{n+1}\right) \right) \\
&&\delta _{1}H\left( T\left( A_{n+1}\right) ,T\left( A_{n+1}\right) \right)
+\delta _{2}H\left( T\left( A_{n+2}\right) ,T\left( A_{n}\right) \right) \\
&\leq &(\alpha +\beta )H\left( T\left( A_{n}\right) ,T\left( A_{n+1}\right)
\right) +\gamma H\left( T\left( A_{n+1}\right) ,T\left( A_{n+2}\right)
\right) , \\
&&\delta _{2}[H\left( T\left( A_{n+2}\right) ,T\left( A_{n+1}\right) \right)
+H\left( T\left( A_{n+1}\right) ,T\left( A_{n}\right) \right) ] \\
&=&\left( \alpha +\beta +\delta _{2}\right) H\left( T\left( A_{n}\right)
,T\left( A_{n+1}\right) \right) +\left( \gamma +\delta _{2}\right) H\left(
T\left( A_{n+1}\right) ,T\left( A_{n+2}\right) \right) .
\end{eqnarray*}%
Now, if $H\left( T\left( A_{n}\right) ,T\left( A_{n+1}\right) \right) \leq
H\left( T\left( A_{n+1}\right) ,T\left( A_{n+2}\right) \right) $, we have%
\begin{eqnarray*}
H(T\left( A_{n+1}\right) ,T\left( A_{n+2}\right) &\leq &\psi \left( \max
\{H\left( T\left( A_{n}\right) ,T\left( A_{n+1}\right) \right) ,H\left(
T\left( A_{n+1}\right) ,T\left( A_{n+2}\right) \right) \}\right) \\
&=&\psi (H\left( T\left( A_{n}\right) ,T\left( A_{n+1}\right) \right) )
\end{eqnarray*}%
for all $n\in 
%TCIMACRO{\U{2115} }%
%BeginExpansion
\mathbb{N}
%EndExpansion
.$ Therefore for $i=1,2,...,n$, we have%
\begin{eqnarray*}
H(T\left( A_{i-1}\right) ,T\left( A_{i}\right) ) &\leq &\psi
(H(A_{i-1},A_{i})), \\
H(T\left( A_{i-2}\right) ,T\left( A_{i-1}\right) ) &\leq &\psi
(H(A_{i-2},A_{i-1})), \\
&&\cdot \cdot \cdot , \\
H(T\left( A_{0}\right) ,T\left( A_{1}\right) ) &\leq &\psi (H(A_{0},A_{1})),
\end{eqnarray*}%
and so we obtain%
\begin{equation*}
H(T\left( A_{n}\right) ,T\left( A_{n+1}\right) \leq \psi
^{n}(H(A_{0},T\left( A_{1}\right) ))
\end{equation*}%
for all $n\in 
%TCIMACRO{\U{2115} }%
%BeginExpansion
\mathbb{N}
%EndExpansion
.$ Follows the similar argument to those in the proof of Theorem 2.1, we get 
$H\left( T\left( A_{n}\right) ,T\left( A_{m}\right) \right) \rightarrow 0$
as $n,m\rightarrow \infty $. Therefore $\{T\left( A_{n}\right) \}$ is a
Cauchy sequence in $T\left( X\right) .$ Since $(T\left( X\right) ,d)$ is
complete in $CB\left( X\right) $, we have $T\left( A_{n}\right) \rightarrow
V $ as $n\rightarrow \infty $\ for some $V\in CB\left( X\right) .$ Also, we
can find $U$ in $CB\left( X\right) $ such that $T(U)=V.$

We claim that $S(U)=T(U).$ If not, then since $\left( T\left( A_{n+1}\right)
,T\left( A_{n}\right) \right) \subseteq E\left( G\right) $ so by property (P$%
^{\ast }$), there exists a subsequence $\{T\left( A_{n_{k}+1}\right) \}$ of $%
\{T\left( A_{n+1}\right) \}$ such that $\left( T\left( U\right) ,T\left(
A_{n_{k}+1}\right) \right) \subseteq E\left( G\right) $ for every $n\in 
%TCIMACRO{\U{2115} }%
%BeginExpansion
\mathbb{N}
%EndExpansion
$. As $\left( U,T\left( U\right) \right) \subseteq E\left( G\right) $ and $%
\left( T\left( A_{n_{k}+1}\right) ,A_{n_{k}}\right) =\left( S\left(
A_{n_{k}}\right) ,A_{n_{k}}\right) \subseteq E\left( G\right) $ implies that 
$\left( U,A_{n_{k}}\right) \subseteq E\left( G\right) .$ Now%
\begin{eqnarray}
H(S\left( U\right) ,T\left( A_{n_{k}+1}\right) ) &=&H(S\left( U\right)
,S\left( A_{n_{k}}\right) )  \notag \\
&\leq &\psi \left( M_{2}\left( U,A_{n_{k}}\right) \right) ,  \TCItag{2.2}
\end{eqnarray}%
where%
\begin{eqnarray*}
M_{2}\left( U,A_{n_{k}}\right) &=&\alpha H(T\left( U\right) ,T\left(
A_{n_{k}}\right) )+\beta H(S\left( U\right) ,T\left( U\right) )+\gamma
H(S\left( A_{n_{k}}\right) ,T\left( A_{n_{k}}\right) ) \\
&&+\delta _{1}H(S\left( U\right) ,T\left( A_{n_{k}}\right) )+\delta
_{2}H(S\left( A_{n_{k}}\right) ,T\left( U\right) ) \\
&=&\alpha H(T\left( U\right) ,T\left( A_{n_{k}}\right) )+\beta H(S\left(
U\right) +T\left( U\right) )+\gamma H(T\left( A_{n_{k}+1}\right) ,T\left(
A_{n_{k}}\right) ) \\
&&+\delta _{1}H(S\left( U\right) ,T\left( A_{n_{k}}\right) )+\delta
_{2}H(T\left( A_{n_{k}+1}\right) ,T\left( U\right) ).
\end{eqnarray*}%
On taking limit as $k\rightarrow \infty $ in (2.2), we have%
\begin{eqnarray*}
H(S\left( U\right) ,T\left( U\right) ) &\leq &\psi \left( \left( \beta
+\delta _{1}\right) H\left( T\left( U\right) ,T\left( U\right) \right)
\right) \\
&<&H(S\left( U\right) ,T\left( U\right) ),
\end{eqnarray*}%
a contradiction. Hence $S\left( U\right) =T\left( U\right) ,$ that is, $U\in
CP\left( S,T\right) $.

\noindent To prove (ii), suppose that $CP\left( S,T\right) $ is complete set
in $G$. Let $U,V\in CP\left( S,T\right) $ and suppose that the
Pompeiu-Hausdorff weight assign to the $S\left( U\right) $ and $S\left(
V\right) $ is not zero. Since pair $(S,T)$ is a graph $\psi _{2}$%
-contraction, we obtain that%
\begin{equation}
H(S\left( U\right) ,S\left( V\right) )\leq \psi (M_{2}(U,V)),  \tag{2.3}
\end{equation}%
where%
\begin{eqnarray*}
M_{2}(U,V)) &=&\alpha H(T\left( U\right) ,T\left( V\right) )+\beta H(S\left(
U\right) ,T\left( U\right) )+\gamma H(S\left( V\right) ,T\left( V\right) ) \\
&&\delta _{1}H(S\left( U\right) ,T\left( V\right) )+\delta _{2}H(S\left(
V\right) ,T\left( U\right) ) \\
&=&\alpha H\left( S\left( U\right) ,S\left( V\right) \right) +\beta
H(S\left( U\right) ,S\left( U\right) )+\gamma H(S\left( V\right) ,T\left(
V\right) ) \\
&=&\left( \alpha +\delta _{1}+\delta _{2}\right) H(S\left( U\right) ,S\left(
V\right) ),
\end{eqnarray*}%
thus%
\begin{eqnarray*}
H(S\left( U\right) ,S\left( V\right) ) &\leq &\psi (\left( \alpha +\delta
_{1}+\delta _{2}\right) H(S\left( U\right) ,S\left( V\right) )) \\
&<&\psi \left( H(S\left( U\right) ,S\left( V\right) )\right) ,
\end{eqnarray*}%
a contradiction as $\psi \left( t\right) <t$ for all $t>0$. Hence (ii) is
proved.

\noindent To prove (iii),\ suppose the set $CP\left( S,T\right) \ $is weakly
compatible.\ First we are to show that $Fix\left( T\right) \cap Fix(S)\ $is
nonempty. Let $W=S\left( U\right) =T\left( U\right) ,$ then we have $T\left(
W\right) =TS\left( U\right) =ST\left( U\right) =S\left( W\right) ,$ which
shows that $W\in CP\left( S,T\right) .$ Thus the Pompeiu-Hausdorff weight
assign to the $S\left( U\right) $ and $S\left( W\right) $ is zero$\ $(by
ii). Hence $W=S\left( W\right) =T\left( W\right) ,$ that is, $W\in Fix\left(
S\right) \cap Fix\left( T\right) .$ Since $CP\left( S,T\right) $ is
singleton set, implies $Fix\left( S\right) \cap Fix\left( T\right) $ is
singleton.

\noindent Finally to prove (iv),\ suppose the set $Fix\left( S\right) \cap
Fix\left( T\right) \ $is complete.\ We are to show that $Fix\left( T\right)
\cap Fix(S)\ $is singleton. Assume on contrary that there exist $U$,$V\in
CB\left( X\right) $ such that $U,V\in Fix\left( S\right) \cap Fix\left(
T\right) $ and $U\neq V$. By completeness of $Fix\left( S\right) \cap
Fix\left( T\right) $,\ there exists an edge between $U$ and $V.$ As pair $%
(S,T)$ is a graph $\psi _{2}$-contraction, so we have%
\begin{eqnarray*}
H(U,V) &=&H(S\left( U\right) ,S\left( V\right) ) \\
&\leq &\psi (M_{2}(U,V)) \\
&=&\psi (\alpha H(T\left( U\right) ,T\left( V\right) )+\beta H(S\left(
U\right) ,T\left( U\right) )+\gamma H(S\left( V\right) ,T\left( V\right) ) \\
&&+\delta _{1}H(S\left( U\right) ,T\left( V\right) )+\delta _{2}H(S\left(
V\right) ,T\left( U\right) )) \\
&=&\psi (\alpha H(U,V)+\beta H(U,U)+\gamma H(V,V)+\delta _{1}H(U,V)+\delta
_{2}H(V,U)) \\
&\leq &\psi \left( H\left( U,V\right) \right) ,
\end{eqnarray*}%
a contradiction.\ Hence $U=V$. Conversely, if $Fix(S)\cap Fix(T)$ is
singleton, then since $E(G)\supseteq \Delta ,$ so it is obvious that $%
F(S)\cap F(T)$ is complete set. $\square $\medskip

\noindent \textbf{Example 2.7.} \ \ Let $X=%
%TCIMACRO{\U{211d} }%
%BeginExpansion
\mathbb{R}
%EndExpansion
_{+}=V\left( G\right) $ be endowed with Euclidean metric $d.$ Let $%
f:X\rightarrow X$ be defined as $f\left( x\right) =\left\{ 
\begin{array}{l}
10,\text{ if }x\in \left[ 0,10\right]  \\ 
20,\text{ otherwise}%
\end{array}%
\right. $ and $\left( a,b\right) \in E\left( G\right) $ for some $a\in A,$ $%
b\in B$ if $b=f\left( a\right) .$ Define $S,T:CB\left( X\right) \rightarrow
CB(X)$ as follows:%
\begin{eqnarray*}
S(U) &=&\left\{ 
\begin{array}{cc}
\left[ 0,10\right] , & if\text{ }U\subseteq \lbrack 0,10] \\ 
\lbrack 10,20], & otherwise%
\end{array}%
\right. \text{ and} \\
T(U) &=&\left\{ 
\begin{array}{cc}
\left[ 0,10\right] , & if\text{ }U\subseteq \lbrack 0,10] \\ 
\lbrack 5,25], & otherwise.%
\end{array}%
\right. 
\end{eqnarray*}%
Note that, for all $V\in CB(X),$ $(S\left( V\right) ,V)\subseteq E\left(
G\right) $ and $(V,T\left( V\right) )\subseteq E\left( G\right) $.

\noindent Let $\psi :%
%TCIMACRO{\U{211d} }%
%BeginExpansion
\mathbb{R}
%EndExpansion
_{+}\rightarrow 
%TCIMACRO{\U{211d} }%
%BeginExpansion
\mathbb{R}
%EndExpansion
_{+}$ be defined by%
\begin{equation*}
\psi \left( \alpha \right) =\left\{ 
\begin{array}{cc}
\dfrac{3}{4}t & 0\leq t<1 \\ 
&  \\ 
\dfrac{5}{6}t, & 1\leq t.%
\end{array}%
\right. 
\end{equation*}%
It is easy to verify that $\psi \in \Psi .$\ Now for all $A,B\in CB\left(
X\right) $ with $S\left( A\right) \neq S\left( B\right) ,$ we consider $%
A\subseteq \lbrack 0,10]$ and $B\nsubseteq \lbrack 0,10]$ with $\left(
A,B\right) \subseteq E\left( G\right) ,$ implies%
\begin{eqnarray*}
H(S\left( A\right) ,S\left( B\right) ) &=&H\left( \left[ 0,10\right]
,[10,20]\right)  \\
&=&10 \\
&<&\frac{100}{9} \\
&=&\psi \left( 15\alpha +5\beta \right)  \\
&=&\psi \left( \alpha H\left( \left[ 0,10\right] ,\left[ 5,25\right] \right)
+\gamma H\left( \left[ 10,20\right] ,\left[ 5,25\right] \right) \right)  \\
&=&\psi \left( \alpha H(T\left( A\right) ,T\left( B\right) )+\gamma H\left(
S\left( B\right) ,T\left( B\right) \right) \right) \leq \psi (M_{2}(A,B)),
\end{eqnarray*}%
where $\alpha =\frac{5}{6},\gamma =\frac{1}{6},$ $\beta =\delta _{1}=\delta
_{2}=0\ $and%
\begin{eqnarray*}
M_{2}(A,B) &=&\alpha H(T\left( A\right) ,T\left( B\right) )+\beta H(S\left(
A\right) ,T\left( A\right) )+\gamma H(S\left( B\right) ,T\left( B\right) ) \\
&&+\delta _{1}H(S\left( A\right) ,T\left( B\right) )+\delta _{2}H(S\left(
B\right) ,T\left( A\right) )
\end{eqnarray*}

\noindent Hence pair $\left( S,T\right) $\ is graph $\psi _{2}$-contraction.
Thus all the conditions of Theorem 2.6 are satisfied. Moreover, the set $%
[0,10]$ is the common fixed point of $S$ and $T$, and $Fix\left( S\right)
\cap Fix\left( T\right) $ is complete. $\square $\medskip 

The following corollary generalizes and extends Theorem 2.1 of \cite{AAKN15}.

\noindent \textbf{Corollary 2.7.} \ \ Let $(X,d)$ be a complete metric space
endowed with a directed graph $G$ such that $V(G)=X$ and $E(G)\supseteq
\Delta .$ Suppose that the mappings $S,T:CB\left( X\right) \rightarrow
CB\left( X\right) $ satisfies the following:

\begin{itemize}
\item[(a)] for every $V\ $in $CB(X),$ $\left( S\left( V\right) ,V\right)
\subset E\left( G\right) $ and $\left( V,T\left( V\right) \right) \subseteq
E\left( G\right) .$

\item[(b)] There exists $\psi \in \Psi $ such that for all $A,B\in CB\left(
X\right) $ with there is an edge between $A$ and $B$ implies%
\begin{equation*}
H(S\left( A\right) ,S\left( B\right) )\leq \psi (\alpha H(T\left( A\right)
,T\left( B\right) )+\beta H(S\left( A\right) ,T\left( A\right) )+\gamma
H(S\left( B\right) ,T\left( B\right) ))
\end{equation*}%
hold, where $\alpha ,\beta ,\gamma $ are nonnegative real numbers with $%
\alpha +\beta +\gamma \leq 1.$ If the range of $T$ contains the range of $S$%
, then the following statements hold:
\end{itemize}

\begin{enumerate}
\item[(i)] $CP\left( S,T\right) \neq \emptyset $ provided that $G$ is weakly
connected with satisfies the property\ (P$^{\ast }$) and $T\left( X\right) $
is complete subspace of $CB\left( X\right) $.

\item[(ii)] if $CP\left( S,T\right) $ is complete, then the
Pompeiu-Hausdorff weight assigned to the $S\left( U\right) $ and $S\left(
V\right) $ is $0$ for all $U,V\in CP\left( S,T\right) $.

\item[(iii)] if $CP\left( S,T\right) $ is complete and $S$ and $T$ are
weakly compatible, then $Fix\left( S\right) \cap Fix\left( T\right) $ is a
singleton.

\item[(iv)] $Fix\left( S\right) \cap Fix\left( T\right) $ is complete if and
only if $Fix\left( S\right) \cap Fix\left( T\right) $ is a singleton.\medskip
\end{enumerate}

\noindent \textbf{Corollary 2.8.} \ \ Let $(X,d)$ be a complete metric space
endowed with a directed graph $G$ such that $V(G)=X$ and $E(G)\supseteq
\Delta .$ Suppose that the mappings $S:CB\left( X\right) \rightarrow
CB\left( X\right) $ satisfies the following:

\begin{itemize}
\item[(a)] for every $V\ $in $CB(X),$ $\left( S\left( V\right) ,V\right)
\subset E\left( G\right) .$

\item[(b)] There exists $\psi \in \Psi $ such that for all $A,B\in CB\left(
X\right) $ with there is an edge between $A$ and $B$ implies%
\begin{equation*}
H(S\left( A\right) ,S\left( B\right) )\leq \psi (\alpha H(A,B)+\beta
H(S\left( A\right) ,A)+\gamma H(B,S\left( B\right) ))
\end{equation*}%
hold, where $\alpha ,\beta ,\gamma $ are nonnegative real numbers with $%
\alpha +\beta +\gamma \leq 1.$ Then the following statements hold:
\end{itemize}

\begin{enumerate}
\item[(i)] if $Fix\left( S\right) $ is complete, then the Pompeiu-Hausdorff
weight assigned to the $U,V\in Fix\left( S\right) $ is $0$.

\item[(ii)] If the weakly connected graph $G$ satisfies the property (P$%
^{\ast }$), then $S$ has a fixed point.

\item[(iii)] $Fix\left( S\right) \ $is complete if and only if $Fix\left(
S\right) \ $is a singleton.
\end{enumerate}

\noindent \textit{Proof.} \ \ If we take $T=I$ (identity map) in above
Corollary 2, the result follows.\medskip

\noindent \textbf{Remark 2.9.}

\begin{itemize}
\item[(1)] If $E(G):=X\times X$, then clearly $G$ is connected and our
Theorem 2.1 improves and generalizes Theorem 2.1 in \cite{AAKN15}, Theorem
2.1 in \cite{BegButt}, Theorem 3.1 in \cite{Jachymski}.

\item[(2)] If $E(G):=X\times X$, then clearly $G$ is connected and our
Theorem 2.4 extends and generalizes Theorem 2.5 in \cite{BegButt}, Theorem
3.2 in \cite{Nadler}, Theorem 5.1 in \cite{Edelstein} and Theorem 3.1 in 
\cite{Jachymski}.

\item[(3)] If $E(G):=X\times X$, then clearly $G$ is connected and our
Corollary 2.5 improves and generalizes Theorem 2.1 in \cite{BegButt},
Theorem 3.2 in \cite{Nadler} and Theorem 3.1 in \cite{Jachymski}.\medskip
\end{itemize}

\noindent \textbf{Conclusion}. Jachymski and Jozwik initiated the study of
ordered structured metric fixed point theory by using the ordered structured
with a graph structure on a metric space. Recently many results appeared in
the literature giving the fixed point problems of mappings endow with a
graph. We presented the common fixed points of a class of multivalued maps
with set-valued domain that are commuting only at their coincidence points
endow with a directed graph. We presented some examples to show the
validated of obtained results.\newline

\end{document}